\documentclass[journal ]{new-aiaa}
\usepackage[utf8]{inputenc}
\usepackage{textcomp}
\usepackage{setspace}
\usepackage[inkscapelatex=false]{svg}

\usepackage{graphicx}
\usepackage{amsmath}
\usepackage[version=4]{mhchem}
\usepackage{siunitx}
\usepackage{longtable,tabularx}
\usepackage{placeins}
\usepackage{makecell}

\usepackage{booktabs}
\usepackage{subcaption}
\usepackage{array}
\usepackage{multirow}
\title{A Semi-analytic Method for Rapid Coverage Analysis of a Satellite Using Piecewise Ellipse Models}

\author[1,*]{Byeong-Un Jo}
\author[2,*]{Jinsung Lee}

\affil[1]{Department of Aerospace System Engineering, Sejong University, Seoul, 05006, Republic of Korea}
\affil[2]{Department of Astronomy and Space Science, Kyung Hee University, Yongin-si, 17104, Republic of Korea}

\newcommand{\corrfootnote}{%
\begingroup
\renewcommand{\thefootnote}{*}
\footnotetext{Corresponding authors, E-mail: bjo@sejong.ac.kr (B.-U., Jo), jinsung.lee@khu.ac.kr (J., Lee)}
\endgroup
}

\begin{document}
\setstretch{1.2}
\maketitle
\corrfootnote

\begin{abstract}
Satellite coverage analysis commonly relies on intensive numerical computations to evaluate access and revisit times, resulting in high computational costs, particularly in long-term performance assessment and large-scale constellation analysis involving numerous satellites. To address this issue, this paper presents a semi-analytic method for rapid coverage analysis for a satellite using piecewise ellipse models. We introduce the epoch longitude of the ascending node (ELAN), defined as the longitude of the ascending node at the beginning of each nodal period. Since the shape of the ground track over one nodal period remains invariant, target visibility within the period is uniquely determined by ELAN. By identifying the feasible range of ELAN, the proposed method excludes invisible nodal periods, thereby substantially reducing search space. Within the feasible range, we transform the visibility-boundary relationship into a modified ELAN (MELAN) and approximate it using analytically invertible piecewise ellipse models, enabling direct evaluation of the entry and exit times without repeated time-domain visibility searches. By exploiting the periodic evolution of ELAN, the proposed method efficiently computes access and revisit times over a given time horizon. Numerical case studies over a wide range of orbital configurations and target latitudes show that the proposed method achieves sub-0.1\% error in most cases while significantly reducing computational effort.\\ 
\end{abstract}


\section{Introduction}
The demand for satellite-based services has increased significantly in recent years, driven by the rapid reduction in launch costs and advancements in small satellite technologies, along with expanding applications such as Earth observation, global communications, navigation, and space-based monitoring\cite{jo2022optimal}. As satellite missions become more diverse and performance requirements more stringent, accurate and rapid evaluation of satellite coverage performance has become increasingly important in mission analysis and system design, particularly in constellation design\cite{ho2024space,jo2024simultaneous}. 

Coverage analysis aims to quantify the visibility and accessibility of a ground target by a satellite or constellation over a given time horizon, providing key performance metrics such as access time, revisit time, and coverage continuity. These metrics should be properly considered in the design and evaluation of space missions, as they directly affect mission effectiveness in applications. In particular, coverage performance is closely related to critical system-level decisions, including orbital altitude, inclination, target location, and sensor characteristics. As a result, coverage analysis is not only used for performance evaluation but also serves as a core component in mission design and trade-space exploration\cite{buzzi2019assessment,Han2019HybridOptimization}. Accordingly, the importance of efficient coverage analysis has increased with the growing complexity of modern space systems, where large design space and long-term mission scenarios must be explored\cite{he2022analytical,zhang2022geometric}. For instance, the evaluation of coverage performance over extended periods, such as weeks to years, is essential for assessing the reliability and sustainability of service provision. In such cases, coverage analysis must be repeatedly performed across a wide range of design variables, making computational efficiency a critical requirement. This is particularly relevant in emerging applications involving large-scale constellations and space-based infrastructure, where rapid and accurate evaluation of coverage metrics is necessary for system-level optimization and real-time decision-making.

Existing approaches can be categorized into numerical, semi-analytic, and geometric/map-based methods. Numerical propagation and grid-based approaches are flexible because they can incorporate detailed orbital dynamics, sensor geometry, and mission constraints. Representative examples include improved grid-point strategies with error control and multiresolution acceleration\cite{Song2018GridPoint} and self-adaptive interpolation methods for rapid visibility determination\cite{han2017rapid}. However, since these methods still rely on temporal or spatial discretization, their computational cost typically grows with the analysis duration, the number of targets, or the number of satellites.

Semi-analytic approaches alleviate part of this burden by expressing coverage metrics more directly in terms of orbital parameters or projected footprint geometry. Sengupta et al. developed semi-analytic techniques for terrestrial coverage as functions of orbital elements\cite{sengupta2010satellite}, and Crisp et al. proposed an elliptical projected-footprint approach for revisit-time estimation\cite{Crisp2018Revisit}. These methods are attractive for early design and optimization, but they are often tailored to specific metrics or geometric approximations and do not generally provide explicit expressions over extended analysis horizons. Furthermore, these methods require suitable initial guesses, which can be disadvantageous in constellation optimization involving extensive design space.

Another important line of research reformulates coverage as a geometric problem in reduced-dimensional space. Ulybyshev established an influential analytic framework for discontinuous coverage analysis\cite{Ulybyshev2015Discontinuous}, while route-theory formulations by Razoumny provided efficient analytic treatments of Earth discontinuous coverage and constellation synthesis\cite{Razoumny2016RouteTheory1}. More recently, Han et al. introduced the field-mapping method for single-satellite ground-target coverage\cite{Han2021FieldMapping}; He and Li combined coverage-region and route-theory concepts for repeating-ground-track constellations\cite{he2022analytical}; Zhang et al. extended the framework to relative and constellation field mapping for multi-satellite revisit analysis\cite{zhang2022geometric}; and subsequent studies further expanded 2-D map methods to minimum-observation prediction, elliptical orbits, and elevation-angle-constrained coverage analysis\cite{Bai2022MinimumObservation,Bai2023Elliptical2DMap,Gu2024Elevation2DMap}. These approaches are highly efficient and intuitive, but they still require auxiliary geometric constructions or branchwise intersection calculations.

The need for rapid yet accurate coverage models is further amplified in constellation optimization and mega-constellation studies, where coverage analysis must be embedded within large iterative search procedures. Examples include hybrid visibility optimization for constellations\cite{Han2019HybridOptimization}, semi-analytic/genetic orbit design\cite{savitri2017satellite}, and rapid coverage analysis methods for mega Walker constellations\cite{Gong2021MegaWalker2D}. In these settings, a formulation that can directly map a compact orbital-phase parameter to access timing information is especially valuable because it weakens the dependence of computational cost on analysis duration and repeated design evaluations. Lee et al. improved the efficiency of regional constellation design by reformulating the problem based on repeating ground tracks (RGTs) and binary integer linear programming\cite{lee2020satellite}. While these approaches significantly improve computational efficiency, they often rely on simplifying assumptions that can limit accuracy and are highly dependent on initial guesses or predefined geometric approximations.  

Motivated by these limitations, this study develops a semi-analytic method for rapid satellite coverage analysis using piecewise ellipse models. To this end, we introduce the epoch longitude of the ascending node (ELAN), defined as the longitude of the ascending node at the beginning of each nodal period, as a discrete orbital-phase parameter that remains constant within each period. Since the ground-track shape over one nodal period is invariant, target visibility within a single nodal period can be parameterized by ELAN, which enables the feasible range of the problem to be identified. Within this range, the visibility-boundary relationship is transformed into a modified ELAN (MELAN), and MELAN is approximated by analytically invertible piecewise ellipse models. This formulation enables direct evaluation of the entry/exit times as functions of ELAN without repeated time-domain visibility searches. By exploiting the periodic evolution of ELAN across successive nodal periods, key coverage metrics, including access time and mean and max revisit times, over extended analysis intervals can then be evaluated systematically and efficiently.

The proposed approach offers several notable advantages over conventional methods. First, since the formulation uses an analytically invertible piecewise-ellipse representation after the visible ELAN range is identified, the computational cost does not increase significantly with the length of the analysis period, unlike conventional methods in which the computation time grows proportionally with the time horizon. Second, as the proposed method does not rely on temporal discretization, it is inherently free from time-step dependency, thereby eliminating the trade-off between accuracy and computational cost. As a result, the proposed method achieves an accuracy comparable to that of high-fidelity simulations with sufficiently small time steps. Third, the formulation remains applicable across a wide range of orbital configurations, making it suitable for repeated coverage evaluations in orbit design and performance analysis (in most cases, the errors are below 0.1\%). Finally, when compared with simplified vectorization-based methods widely used in constellation design, the proposed approach achieves a substantial reduction in computational time. Numerical results indicate that the computation time can be reduced by several orders of magnitude, 1/15,000 of that required by the simplified vectorization-based method in certain cases, which represents a significant improvement for large-scale coverage analysis and design applications.

\section{Problem Definition}
\subsection{Satellite-Target Visibility Formulation}
Satellite coverage analysis fundamentally relies on determining visibility of a given ground target from a satellite at a specific time, which is equivalent to assessing whether the target lies within the field of view of the satellite-mounted sensor. This section introduces the satellite-target visibility problem for which the proposed semi-analytic formulation is developed.

The motion of a satellite is described by the classical six orbital elements: semi-major axis ($a$), inclination angle ($i$), eccentricity ($e$), argument of perigee ($\omega$), right ascension of the ascending node (RAAN, $\Omega$), and true anomaly ($\nu$). While the semi-major axis, eccentricity, and inclination angle are treated as constant in the mean sense over the time span of interest, the RAAN and the argument of perigee exhibit a secular variation caused by the Earth’s oblateness. In the context of coverage analysis, the long-term evolution of these elements is characterized using their mean orbital dynamics\cite{vallado2001fundamentals}, expressed as
\begin{align}
    \dot{\omega} &=
    \frac{3}{2} J_2 \left( \frac{R_\oplus}{p} \right)^2
    \sqrt{\frac{\mu_\oplus}{a^3}}
    \left( 2 - \frac{5}{2} \sin^2 i \right) \label{wdot}
    \\
    \dot{\Omega} &=
    -\frac{3}{2} J_2 \left( \frac{R_\oplus}{p} \right)^2
    \sqrt{\frac{\mu_\oplus}{a^3}} \cos i
    \\
    \dot{M} &=
    \sqrt{\frac{\mu_\oplus}{a^3}}
    \left[1 - \frac{3}{2} J_2 \left( \frac{R_\oplus}{p} \right)^2
    \sqrt{1 - e^2}
    \left( \frac{3}{2} \sin^2 i - 1 \right)
    \right]
\end{align}
where $R_\oplus$, $\omega_\oplus$, and $\mu_\oplus$ are the radius, angular rate, and the standard gravitational parameter of the Earth, respectively. In addition, $\dot{M}$ is the rate of mean anomaly, $p = a(1 - e^2)$ is the semi-latus rectum, and $J_2$ is the dominant zonal harmonic coefficient that represents the effect of the Earth’s equatorial bulge on its gravitational field.

For target visibility analysis, it is more convenient to transform the orbital elements into an Earth-centered, Earth-fixed (ECEF) ground track and perform the analysis in that domain, rather than directly analyzing the orbital elements. Accordingly, the instantaneous latitude ($\phi$) and longitude ($\lambda$) of the satellite in the ground track are expressed explicitly as follows
\begin{align}
    \sin{\phi} &= \sin{i} \sin{u}, & & \hspace{-10pt} \phi \in [-90^\circ, 90^\circ] \label{eq:lat0} \\[3pt]
    \lambda &= \Omega - \theta_{G} 
    + \tan^{-1} (\cos{i} \tan{u}), & & \hspace{-10pt} \lambda \in [-180^\circ, 180^\circ] \\
    &= \Omega_{0} - \theta_{G_0} + (\dot{\Omega}- \omega_{\oplus})t + \tan^{-1} (\cos{i} \tan{u})
    \label{eq:lon0}
\end{align}
where $u=\nu+\omega$ is the argument of latitude, $\theta_{G}$ is the Greenwich mean sidereal time (GMST), and $t$ is the elapsed time from a reference epoch. Additionally, the subscript 0 denotes the initial value at a reference epoch. For a circular orbit, we can reasonably alternate the true anomaly with the mean anomaly, and therefore, the time derivative of the argument of latitude is as follows
\begin{align}
    \dot{u} = \dot{M}+\dot{\omega} \label{udot}
\end{align}

The visibility between the satellite and the ground target is determined by a geometric constraint that accounts for the Earth's curvature and the sensor's field of view. A visibility function is defined as follows: 
\begin{align}
V(t) &= \cos{\phi}\cos{\phi_{t}}\cos{(\lambda-\lambda_{t})} + \sin{\phi}\sin{\phi_{t}} - \cos{\beta} 
\label{eqV}
\end{align}
where the subscript $t$ indicates the values of the target. In Eq. \eqref{eqV}, $\beta$ is the Earth central angle corresponding to the sensor half-angle ($\gamma_{\mathrm{half}}$), representing the angular radius of the ground coverage region measured from the Earth's center, expressed as
\begin{align}
\beta
= -\gamma_{\mathrm{half}}
+ \sin^{-1}\!\left( \frac{r}{R_\oplus}\,\sin\gamma_{\mathrm{half}} \right) \label{eq:beta}
\end{align}
In Eq. \eqref{eq:beta}, $r$ denotes the geocentric radius of the satellite. The target is visible when the function is non-negative; otherwise, it is invisible.
\begin{equation}
    \begin{cases} 
      \text{Visible}, & \text{if } V(t) \geq 0, \\
      \text{Invisible}, & \text{otherwise.}
    \end{cases}
\end{equation}
\\
The instants at which this visibility function becomes zero correspond to the entry and exit times to the visible region. Accordingly, by solving this root-finding problem with the given dynamics of a satellite in Eqs. \eqref{wdot}-\eqref{udot}, we can compute the entry/exit times to the visible region. However, this procedure is often computationally expensive, particularly for large satellite constellations or extended analysis horizons, and its convergence characteristics can be sensitive to the initial guess.

\subsection{Introduction to Epoch Longitude of Ascending Node, ELAN($\lambda_{\Omega}$)}

\begin{figure}[hbt!]
\centering
\includegraphics[width=0.5\columnwidth]{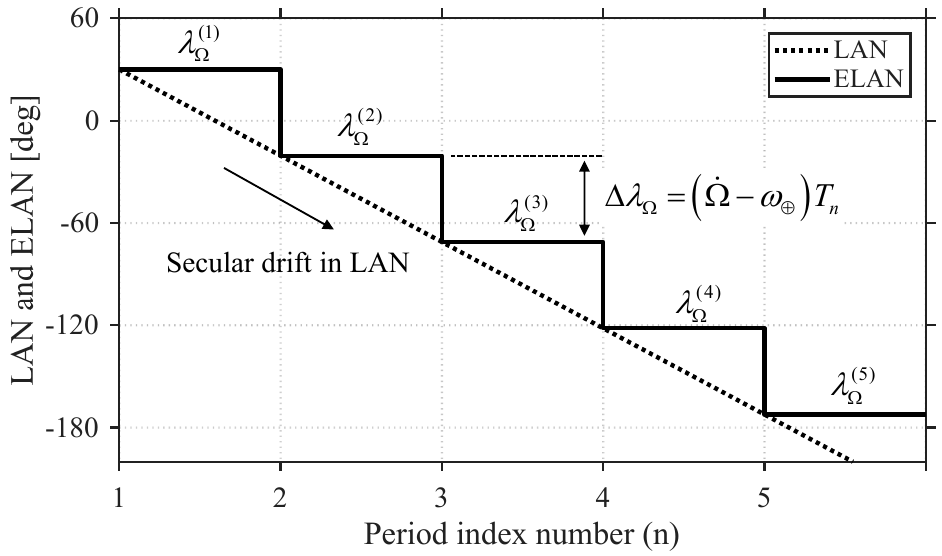}
\caption{Comparison between LAN and ELAN}
\label{Fig:LAN&ELAN}
\end{figure}
This section defines Epoch Longitude of the Ascending Node (ELAN), denoted by $\lambda_{\Omega}$ as the longitude of the ascending node (LAN) at the beginning of the $n$-th nodal period. Unlike LAN, which evolves over time due to the rotation of the Earth and J2 effect, ELAN remains constant within a single nodal period and is updated at each subsequent node as shown in Fig. \ref{Fig:LAN&ELAN}. Therefore, ELAN is a discrete value updated at every nodal period defined as
\begin{align}
\lambda_{\Omega}^{(n)} &= \Omega_0 - \theta_{G_0} + (\dot{\Omega} - \omega_\oplus)\, T_n\, (n-1),
\quad n = 1, 2, 3, \dots, \left\lfloor \frac{t_f}{T_n} \right\rfloor + 1. \label{eq:ELAN}
\end{align}
where $n$ is the index number of $n$-th nodal period, $t_f$ is the final time of interest, and $T_n$ is the nodal period, the time interval between successive crossings of the ascending node, defined as,
\begin{align}
    T_n = \frac{2\pi}{\dot{\omega}+\dot{M}}
\end{align}

\begin{figure}[hbt!]
\centering
\includegraphics[width=0.5\columnwidth]{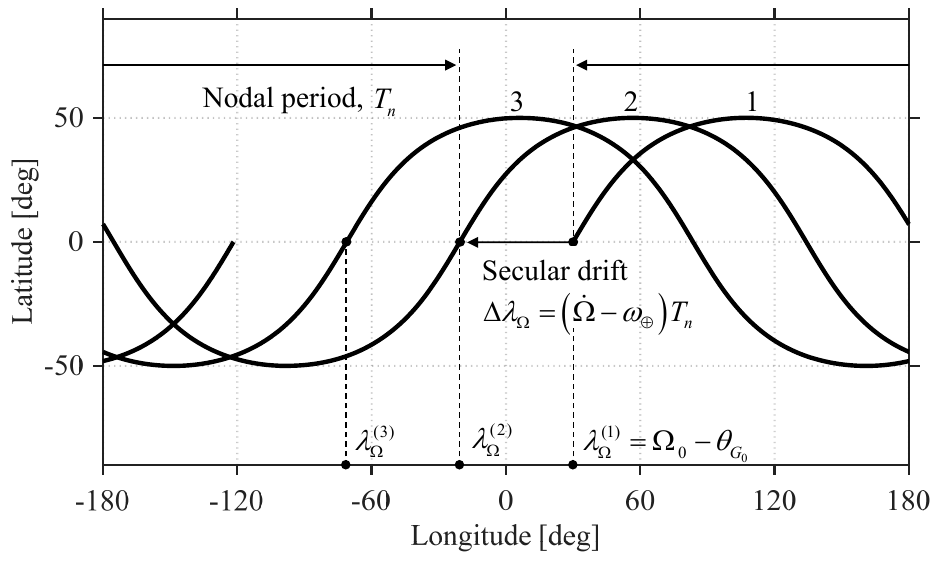}
\caption{Schematic diagram explaining the definition of ELAN}
\label{Fig:ELAN}
\end{figure}
Note that the nodal period is distinct from the commonly used orbital period (Keplerian period). Once the initial RAAN ($\Omega_0$) and GMST ($\theta_{G_0}$) are given, we can simply predict $\lambda_{\Omega}$ for every period as its time derivative is constant. Therefore, if the explicit expressions for the relationship between the entry/exit times to the visible region and ELAN is available, we can directly compute the access time for the given $\lambda_{\Omega}$ from the difference between the entry and exit times. As a result, it is possible to evaluate the total access time over the entire analysis horizon efficiently.

\begin{figure}[hbt!]
\centering
\includegraphics[width=0.6\columnwidth]{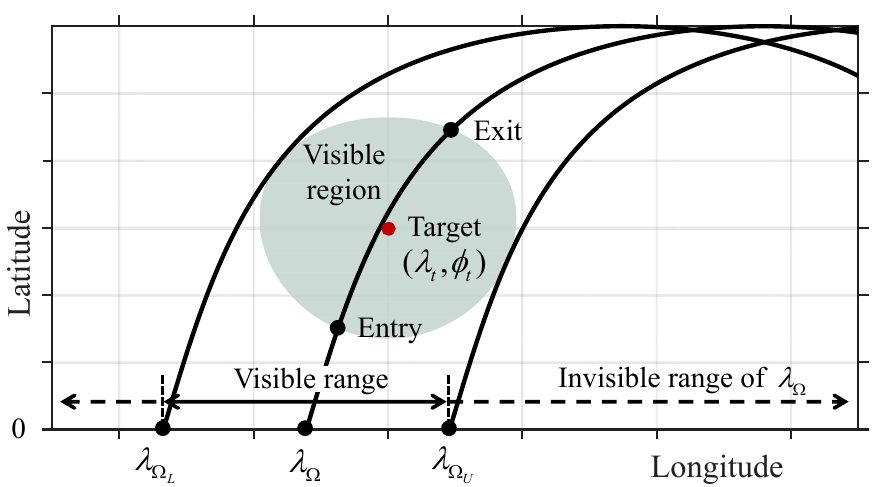}
\caption{Schematic diagram explaining the definition of ELAN and its visible range}
\label{Fig:ELAN0}
\end{figure}
Fig. \ref{Fig:ELAN} explains the definition of $\lambda_{\Omega}$ in the geographic coordinate system. Once ELAN is specified, the ground-track geometry over one nodal period is uniquely determined. Additionally, variations in ELAN result only in a longitudinal shift of the ground track, without altering its geometric shape. This property leads to two important advantages. First, the set of ELANs that can yield target visibility is confined to a limited feasible range, denoted by $[\lambda_{\Omega_L}, \lambda_{\Omega_U}]$ as shown in Fig. \ref{Fig:ELAN0}. For the values outside this range, the corresponding ground tracks never intersect the visible region of the target. Therefore, we can exclude the ground tracks associated with ELANs outside the feasible range from the visibility analysis, significantly reducing the calculation domain.

Second, since the ground-track shape is invariant with respect to ELAN, the visibility assessment can be reformulated based on the geometric relationship between a single representative ground track and the visible region. Once this geometric relationship is characterized, it is no longer necessary to evaluate the visibility condition in Eq. \eqref{eqV} at every time instant along each ground track.

Specifically, if explicit expressions for access time as a function of ELAN is available, we can compute the access time ($t_a$) corresponding to any ELAN within the visible range without performing time-domain visibility checks along the ground track. Accordingly, the access time at the $(n)-th$ nodal period can be expressed as
\begin{equation}
    t_a^{(n)} = 
    \begin{cases} 
      f_a(\lambda_{\Omega}^{(n)}), & \text{if } \space \lambda_{\Omega_L}\leq \lambda_{\Omega}^{(n)} \leq\lambda_{\Omega_U} \\
      0, & \text{otherwise.}
    \end{cases}
\end{equation}
where $f_a(\lambda_{\Omega}^{(n)})$ denote the access-time function with respect to ELAN. Then, using the temporal evolution of ELAN over multiple nodal periods in Eq. \eqref{eq:ELAN}, we can efficiently calculate the total access time over the entire analysis span consisting of the final nodal period $(n_f)$ as $\sum_{n=1}^{n_f} t_a^{(n)}$, enabling rapid evaluation of long-term coverage performance.

\begin{figure}[hbt!]
\centering
\includegraphics[width=0.6\columnwidth]{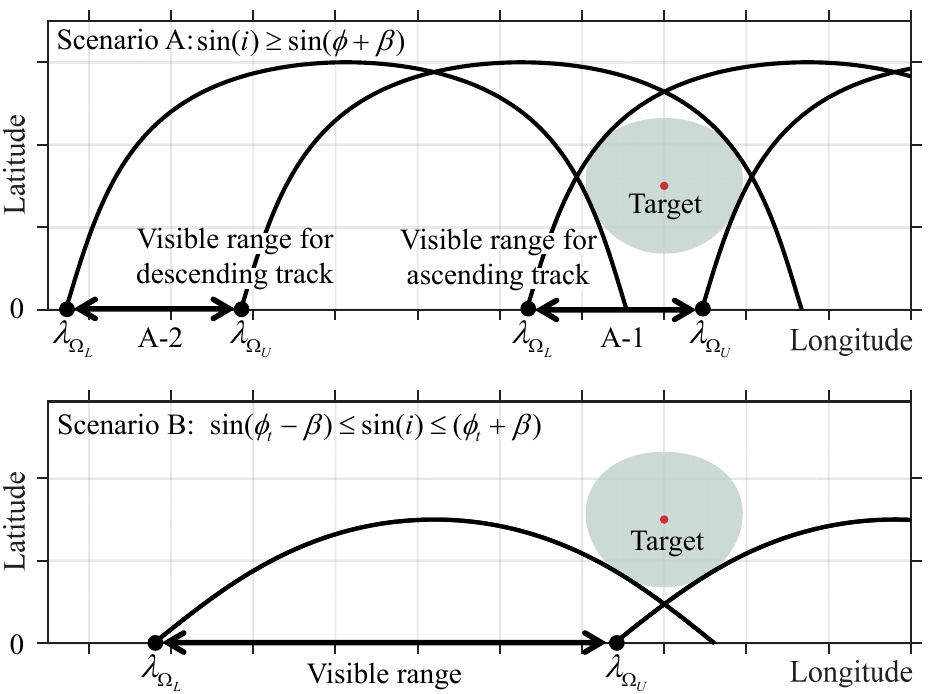}
\caption{Possible Scenarios for the visible range of ELAN}
\label{Fig:Cases}
\end{figure}

Depending on the relationship between the orbital inclination and the target latitude, two main cases exist regarding target accessibility. In Case A, the satellite accesses the target either along the ascending ground track or along the descending ground track, and the access is clearly separable between the two. In contrast, in Case B, access is possible on the ascending track, on the descending track, or on both tracks. Fig. \ref{Fig:Cases} depicts the possible cases for the visible range of ELAN.

Assuming that the ELAN lies within the visible range, when the inclination satisfies $\text{sin}(i) > \text{sin}(\phi_t + \beta)$, the satellite traverses the visible region entirely during both the ascending and descending ground tracks, in which case the entry and exit times are well defined for each track. In contrast, when $\text{sin}(\phi_t - \beta)\le \text{sin}(i) \le \text{sin}(\phi_t + \beta)$, the satellite enters the visible region during the ascending track and exits during the descending track, or the entry and exit times could be defined for only one (descending or ascending) of the two tracks.

Accordingly, three distinct cases are identified: A-1) $\text{sin}(i) > \text{sin}(\phi_t + \beta)$ with access along the ascending track, A-2) $\text{sin}(i) > \text{sin}(\phi_t + \beta)$ with access along the descending track, and B) $\text{sin}(\phi_t - \beta)\le \text{sin}(i) \le \text{sin}(\phi_t + \beta)$. Although the proposed approach remains applicable to Case B, these configurations are generally not preferred in practical mission design due to their limited access characteristics. Therefore, for clarity and practical relevance, this study focuses exclusively on Case A.

\subsection{Problem Definition}
The primary goal of this study is to establish a functional relationship between the access time and ELAN in the form of \( t_a = f_a(\lambda_{\Omega}) \). To compute the access time, the entry and exit times must first be determined, since the access time is obtained from their difference as
\begin{align}
    t_a = t_{\mathrm{out}} - t_{\mathrm{in}} .
\end{align}
Here, \( t_{\mathrm{in}} \) and \( t_{\mathrm{out}} \) denote the times at which the satellite enters and exits the visible region, respectively. These are collectively referred to as visibility-boundary-crossing times. For conciseness, we define the visibility-boundary-crossing time as
\begin{equation}
    t_b \in \{ t_{\mathrm{in}},\ t_{\mathrm{out}} \} .
\end{equation}
Before deriving the functional relationship, we first identify the visible (feasible) range of ELAN, denoted by \( [\lambda_{\Omega_L}, \lambda_{\Omega_U}] \). Let \( \lambda_{\Omega} \) be the independent variable, and consider the problem of determining $t_b$ over a single nodal period. It should be noted that the proposed formulation solves the visibility problem within a single nodal period specified by a given ELAN. Therefore, the quantities \( t_a \) and \( t_b \) defined in this study represent local times defined within a single nodal period over the interval \( [0, T_n] \).

Since the shape of the ground track remains identical for each period, a representative period is selected by setting \( \omega_0 = 0 \), \( t_0 = 0 \), and $\theta_{G_0}=0$. Accordingly, the expressions for the instantaneous latitude and longitude of a satellite during a single period in Eqs. \eqref{eq:lat0}-\eqref{eq:lon0} become
\begin{align}
    u &= \dot{u}t=(\dot{M}+\dot{\omega})t, & & \hspace{-70pt} t \in [0, T_n] \label{eq:u}\\
    \sin{\phi} &= \sin{i} \sin{u}, & & \hspace{-70pt} \phi \in [-90^\circ, 90^\circ] \label{eq:lat} \\
    \lambda &= \Omega_0 + (\dot{\Omega}- \omega_{\oplus})t + \tan^{-1} (\cos{i} \tan{u}), & & \hspace{-70pt} \lambda \in [-180^\circ, 180^\circ] \label{eq:lon} \\[3pt]
    &= \lambda_{\Omega} + \lambda_{l}. \notag & &
\end{align}
where $\lambda_{l}(\triangleq(\dot{\Omega}- \omega_{\oplus})t + \tan^{-1} (\cos{i} \tan{u}))$ denotes the local longitude component relative to the orbital plane. It represents only the longitude variation induced by the orbital motion, J2 effect, and Earth’s rotation.  The initial RAAN ($\Omega_{0}$) in Eq. \eqref{eq:lon} is typically treated as a given parameter. However, since the problem in this subsection is confined to a single nodal period, we replace ($\Omega_{0}$) with $\lambda_{\Omega}$, treating it as a variable. As a satellite crosses the equator at a different longitude in each period, this replacement allows the analysis of how access time varies depending on ELAN ($\lambda_{\Omega}$) at every period. Putting Eq. \eqref{eq:lon} into the visibility function in Eq. \eqref{eqV} yields,
\begin{align}
V(t,\lambda_{\Omega}) &= \cos{\phi}\cos{\phi_{t}}\cos{(\lambda_{\Omega} + \lambda_{l}-\lambda_{t})} + \sin{\phi}\sin{\phi_{t}} - \cos{\beta} \label{eq:NewV}
\end{align}
According to the visibility condition, when Eq. \eqref{eq:NewV} becomes zero, the satellite either enters or exits the visible region. Thus, the visibility-boundary-crossing-time $(t_b)$ can be obtained by finding the root of Eq. \eqref{eq:NewV} for a given $\lambda_{\Omega}$. Although solving this root-finding problem is not inherently challenging, the variations in $\lambda_{\Omega}$ necessitate solving this equation for each nodal period, resulting in a substantial computational load, particularly as the analysis time horizon increases. To mitigate this issue, this study proposes explicit expressions for $t_b (\in \{ t_{\mathrm{in}},\ t_{\mathrm{out}} \})$ as functions of $\lambda_{\Omega}$ using piecewise ellipse models with ten informative data points based on  Eqs. \eqref{eq:u}-\eqref{eq:NewV}. Ultimately, the access time and mean and max revisit times can be computed from $t_{\text{in}}^{(n)}$ and $t_\text{{out}}^{(n)}$.

\section{Determination of Visible Range}
\begin{figure}[hbt!]
\centering
\includegraphics[width=0.7\columnwidth]{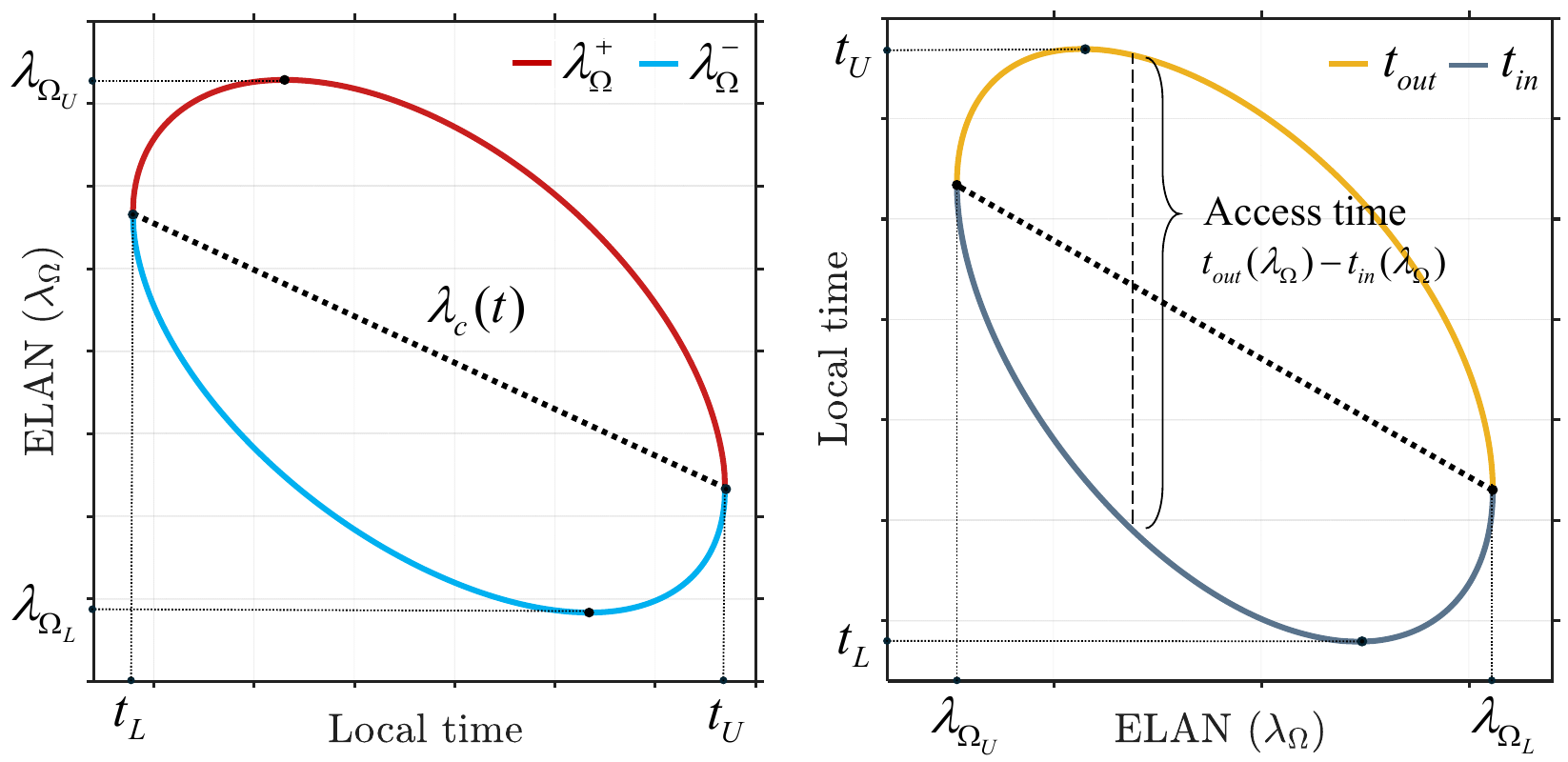}
\caption{Example of ELAN versus visibility-boundary-crossing time and its inverse}
\label{Fig_visible_boundary}
\end{figure}

This section mathematically determines the visible (feasible) range of the key parameters, namely $t_b$ and $\lambda_{\Omega}$, thereby excluding unnecessary analyses in the invisible range. We begin by establishing the relationship between $\lambda_{\Omega}$ and $t_b$ by setting Eq. \eqref{eq:NewV} to zero, as follows
\begin{align}
\lambda_{\Omega}=
\begin{cases}
\lambda_{\Omega}^+(t_{b}) = \lambda_{t} - \lambda_{l} + \cos^{-1}{\dfrac{C_1}{C_2}},\\
\lambda_{\Omega}^-(t_{b}) = \lambda_{t} - \lambda_{l} - \cos^{-1}{\dfrac{C_1}{C_2}},
\label{eq:time2lambda}
\end{cases}
\end{align}
where $C_1 = \cos{\beta} - \sin{\phi}\sin{\phi_t}$ and $C_2 = \cos{\phi_{t}}\cos{\phi}$. Eq. \eqref{eq:time2lambda} defines a visible boundary within which the values satisfy the visibility condition, as illustrated in Fig. \ref{Fig_visible_boundary}. Fig. \ref{Fig_visible_boundary} illustrates an example of both $\lambda_{\Omega}(t_b)$ and its inverse function, $t_b(\lambda_{\Omega})(\in [t_{\text{in}}(\lambda_{\Omega}),\ t_{\text{out}}(\lambda_{\Omega})])$, which are obtained by plotting Eq. \eqref{eq:time2lambda} in local time and ELAN coordinate systems, respectively. In the term with the $\pm$ sign, the region above the straight line connecting the two points for the minimum  ($t_L$) and the maximum ($t_U$) visibility-boundary-crossing time corresponds to the ‘$+$’ sign, while the region below corresponds to the ‘$-$’ sign.

Eq. \eqref{eq:time2lambda} provides an explicit formulation of $\lambda_{\Omega}$ as a function of $t_b$. However, the primary objective of this study is to express $t_b(\lambda_{\Omega})$ explicitly. This could in principle be achieved from Eq. \eqref{eq:time2lambda}, but obtaining its inverse analytically is challenging due to the highly nonlinear nature. Alternatively, this study proposes an alternative approach by reformulating $\lambda_{\Omega}(t_b)$ into a structure that allows for the derivation of its inverse function, such as an elliptical equation. While it is theoretically possible to derive an approximated form of $t_b(\lambda_{\Omega})$ directly by regression, it is more practical in terms of computational cost to establish the invertible expression for $\lambda_{\Omega}(t_b)$ and then determine its inverse in an analytic form. This approach is advantageous in that the analytic forms for the feasible range of time (i.e., its minimum and maximum values) are available, and the maximum value of $\lambda_m$, which is described in the following section, can be determined merely by solving a root-finding problem.

\subsection{Visible range of time}
The bounds of visible time $([t_{L}, t_{U}])$ at which the satellite crosses the visible region are uniquely defined for given $a, e$, and $i$. As depicted in Fig. \ref{Fig_visible_boundary}, within a single nodal period, the satellite is unable to traverse the visible region beyond these bounds, irrespective of $\lambda_{\Omega}$. Therefore, it is necessary to limit the search range by determining the bounds to eliminate infeasible solutions. In Eq. \eqref{eq:time2lambda}, when \(\lvert C_1/C_2 \rvert > 1\), the solution becomes infeasible. Therefore, setting $C_1/C_2 = \pm 1$ provides the bounds as follows:
\begin{align}
& C_1 = \pm C_2, \nonumber\\
=> & \cos{\beta} - \sin{\phi} \sin{\phi_t} = \pm \cos{\phi_t} \cos{\phi}, \nonumber\\
=> & 
    \begin{cases} 
        \phi_{L}=\sin^{-1} \left( \cos{\beta} \sin{\phi_t} 
        - \sqrt{\cos^2{\beta} \sin^2{\phi_t} - \cos^2{\beta} + \cos^2{\phi_t}} \right) = \phi_t - \beta  \\[8pt]
        \phi_{U}=\sin^{-1} \left( \cos{\beta} \sin{\phi_t}  
        + \sqrt{\cos^2{\beta} \sin^2{\phi_t} - \cos^2{\beta} + \cos^2{\phi_t}} \right)  = \phi_t + \beta
    \end{cases} \label{eq:latminmax}.
\end{align}
Note that the maximum latitude for Case B where $\text{sin}(\phi_t - \beta)\le \text{sin}(i) \le \text{sin}(\phi_t + \beta)$ is constrained by the inclination angle as $\phi_{max}=i$. However, rather than explicitly imposing it, it is more convenient to manage the time span directly. Rearranging Eq. \eqref{eq:lat} after incorporating Eq. \eqref{eq:latminmax} gives the lower and upper bounds of the visible time as:
\begin{align}
\text{Case A-1}: \qquad
t_{L} &=\frac{\sin^{-1} \left( \frac{\sin{\phi_{L}}}{\sin{i}} \right)}{\dot{M}+\dot{\omega}},\qquad
t_{U} =\frac{\sin^{-1} \left( \frac{\sin{\phi_{U}}}{\sin{i}} \right)}{\dot{M}+\dot{\omega}} \nonumber \\
\text{Case A-2}: \qquad
t_{L} &=\frac{\pi-\sin^{-1} \left( \frac{\sin{\phi_{U}}}{\sin{i}} \right)}{\dot{M}+\dot{\omega}},\qquad
t_{U} =\frac{\pi-\sin^{-1} \left( \frac{\sin{\phi_{L}}}{\sin{i}} \right)}{\dot{M}+\dot{\omega}} \\
\text{Case B}: \qquad
t_{L} &=\frac{\sin^{-1} \left( \frac{\sin{\phi_{L}}}{\sin{i}} \right)}{\dot{M}+\dot{\omega}},\qquad
t_{U} =\frac{\pi-\sin^{-1} \left( \frac{\sin{\phi_{U}}}{\sin{i}} \right)}{\dot{M}+\dot{\omega}} \nonumber
\end{align}
For Case A, two groups of bounds exist: the accesses along the ascending and descending tracks. The corresponding ELANs ($\lambda_{\Omega}(t_L)$ and $\lambda_{\Omega}(t_U)$) for the bounds can be calculated by putting $t_L$ and $t_U$ into Eqs. \eqref{eq:u}-\eqref{eq:time2lambda}. 

\subsection{Visible range of ELAN ($\lambda_{\Omega}$)}
If a satellite crosses the equator beyond the visible range of $\lambda_{\Omega}$, it cannot observe the target during that period. Therefore, in order to improve computational efficiency, it is essential to establish the feasible range of $\lambda_{\Omega}$, along with the visible range of time. Notably, orbits with $\lambda_{\Omega}$ outside the visible range are not considered for analysis, as they cannot access the target; this significantly simplifies the overall analysis. We can determine the range by solving the following optimization problems, which seek to maximize and minimize $\lambda_{\Omega}$.\\\\
$[$Problem 1$]$ - Optimization problems for the feasible range of $\lambda_{\Omega}^+$ and $\lambda_{\Omega}^-$, respectively,
\begin{align}
    \phantom{= \lambda_{\Omega}, \qquad t \in [t_{L}, t_{U}]} 
    & 1) \max_{t_b} J = \lambda_{\Omega}^+(t_b),\qquad 2) \min_{t_b} J = \lambda_{\Omega}^-(t_b),
    & t_b \in [t_{L}, t_{U}] \label{eq:OptP1}
\end{align}
The search range is confined by the visible time span $([t_{L}, t_{U}])$. We can solve the optimization problems in Eq. \eqref{eq:OptP1} by directly applying optimization techniques such as the golden section search or by differentiating $\lambda_{\Omega}$ with respect to $t_b$, transforming the problem into finding the root of the derivative ($\text{d}\lambda_{\Omega}(t_b)/\text{d}t_b$), and then using root-finding techniques to solve it. Since the first and second derivatives of $\lambda_{\Omega}$ for $t_b$ can be analytically derived, a computationally efficient technique, the Newton-Raphson method is applicable, but it may diverge if the initial guess is far from the true root. Furthermore, since the solution is likely biased toward the boundary, employing a bracketing root-finding method would be inefficient. Given that the bounds exist and $\lambda_{\Omega}(t_b)$ is a unimodal function for $t_b$ within the bounds as illustrated in Fig. \ref{Fig_visible_boundary}, the golden-section search method is a suitable technique in that it guarantees convergence. Consequently, solving Problem 1 provides the bounds of visible ELANs ($\lambda_{\Omega_L}$ and $\lambda_{\Omega_U}$) and the corresponding times ($t_{\Omega_L}$ and $t_{\Omega_U}$).

\subsection{Modified ELAN (MELAN) and its visible range}
Since Eq. \eqref{eq:time2lambda} is highly nonlinear and does not have a form suitable for direct analytic inversion, converting it into an analytically invertible function without specific manipulation requires a large number of data points. To address this issue, this subsection modifies $\lambda_{\Omega}(t_b)$ to an invertible modified variable named the Modified Epoch Longitude of the Ascending Node (MELAN), by subtracting the straight line connecting two points for $t_L$ and $t_U$. Let $\lambda_{c}$ represent the line connecting two points for $t_L$ and $t_U$ in Fig. \ref{Fig_visible_boundary} expressed as
\begin{align}
\lambda_{c}(t_b) = k_1t_b+k_2=\frac{\lambda_{\Omega}(t_U)-\lambda_{\Omega}(t_L)}{t_U-t_L}t_b + \frac{t_U\lambda_{\Omega}(t_L)-t_L\lambda_{\Omega}(t_U)}{t_U-t_L} \label{eq:lambda_c}
\end{align}
Subtracting $\lambda_c$ from $\lambda_{\Omega}$ yields the MELAN, $\lambda_m$, as follows
\begin{align}
\lambda_{m}=
\begin{cases}
\lambda_{m}^+ = \lambda_{\Omega}^+ - \lambda_c\\
\lambda_{m}^- = \lambda_{\Omega}^- - \lambda_c    
\end{cases}
\label{eq:MELAN}
\end{align}

MELAN, $(\lambda_m)$, exhibits strong elliptical characteristics. Although the modified function does not constitute a perfect ellipse, its elliptical characteristics enable regression with a significantly reduced number of data points while ensuring acceptable accuracy. 
\begin{figure}[hbt!]
\centering
\includegraphics[width=1\columnwidth]{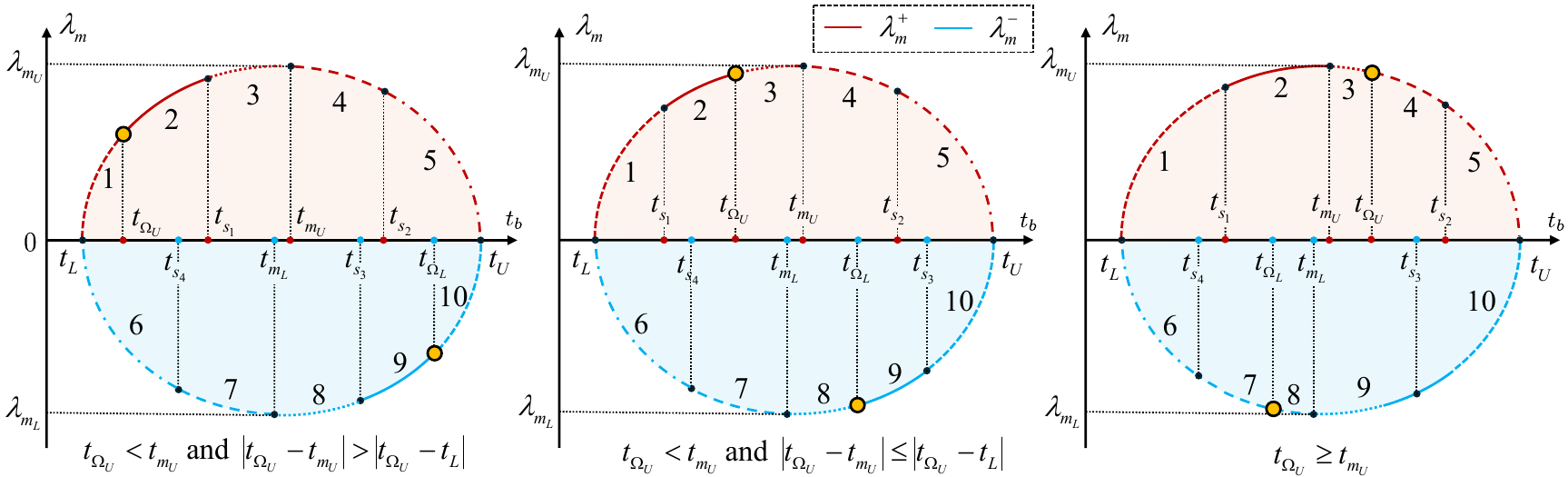}
\caption{Examples of $\lambda_{m}$: three possible cases} 
\label{Fig:MELAN}
\end{figure}
Fig. \ref{Fig:MELAN} illustrates an example of $\lambda_m$ with the corresponding index number of $i$-th model described in the following section. Analyzing Eq. \eqref{eq:MELAN} and Fig. \ref{Fig:MELAN} reveals that $\lambda_m$ exhibits strong elliptical characteristics for the following reasons:\\
\hspace*{0.3cm} 1) Major and minor axes: the function forms an elongated shape with clearly distinguishable major and minor axes,\\
\hspace*{0.3cm} 2) Tangent properties at extremes: the slope is zero or undefined at the vertex and co-vertex.\\

While a substantial number of data points can be extracted from Eq. \eqref{eq:MELAN} with respect to various $t_b$ for a perfect fit, it is possible to achieve an acceptable level of accuracy with significantly fewer data points by leveraging the elliptical characteristics. Since $\lambda_{m}$ does not form a perfect ellipse, a more precise fitting requires separate piecewise ellipse models. Each model is distinguished by ten data points: \\
\hspace*{0.3cm} 1) two points corresponding to the min/max time $[(t_L,0),(t_U,0)]$,\\
\hspace*{0.3cm} 2) two points corresponding to the min/max ELAN $[(t_{\Omega_U},\lambda_m(t_{\Omega_U}),(t_{\Omega_L},\lambda_m(t_{\Omega_L})]$\\
\hspace*{0.3cm} 3) two points for min/max MELAN $[(t_{m_L},\lambda_{m_L}), (t_{m_U},\lambda_{m_U})]$ , which constitute the visible range of $\lambda_m$. \\
\hspace*{0.3cm} 4) four points with supplementary data for gap compensation\\
The first six points have mathematically important characteristics in forming a piecewise ellipse model as they correspond to the visible limits. Among them, the first four points are determined in the previous sections, and this section describes an optimization approach to find the visible range of MELAN (min/max values of MELAN) for the next two points. Accordingly, we define the following optimization problems
\begin{align}
    \phantom{= \lambda_{\Omega}, \qquad t \in [t_{L}, t_{U}]} 
    & 1) \max_{t_b} J = \lambda_{m}^+(t_b),\qquad 2) \min_{t_b} J = \lambda_{m}^-(t_b),
    & t_b \in [t_{L}, t_{U}] \label{eq:OptProblem}
\end{align}
The objective functions are to maximize and minimize $\lambda_{m}^+$ and $\lambda_{m}^-$, respectively, within the search range constrained by the visible time span $([t_{L}, t_{U}])$. In a similar manner to solving Problem 1, we can solve the optimization problems in Eq. \eqref{eq:OptProblem} by applying the golden-section search. However, unlike Problem 1, these optimization problems benefit from a fine initial guess as the min/max values occur near the midpoint in the visible range of time. Moreover, the analytic expressions of the time derivative are available. Therefore, instead of directly finding the extreme values, we find the root of its derivative by adopting the Newton-Raphson method for computational efficiency with a fine initial guess.\\\\
$[$Problem 2$]$ Root-finding problems for the visible range (min/max range) of $\lambda_{m}$:
\begin{align}
    \text{find $t^*_b$ such that} \,\,\frac{\text{d}\lambda_{m}(t_{b}^*)}{\text{d}t_b} = 0, \quad \text{for } t_b \in [t_{L}, t_{U}]. \label{eq:RFProblem}
\end{align}
In Eq. \eqref{eq:RFProblem}, $\text{d}\lambda_{m}/\text{d}t_b$ is obtained by differentiating Eq. \eqref{eq:MELAN} for $t_b$ as follows:
\begin{align}
&\dfrac{\text{d}\lambda_{m}(t_{b})}{\text{d}t_b} = \dot{\lambda}_{m} = \dot{\lambda}_{\Omega}-\dot{\lambda}_{c}\ = \left( \frac{t_{\phi_t}}{c_{\phi}s_{\lambda_{\Delta}}} - \frac{c_{\beta}t_{\phi}}{c_{\phi_t}c_{\phi}s_{\lambda_{\Delta}}} \right) \dot{\phi} - \dot{\lambda}_l - \dot{\lambda}_c \label{eq:dlam}
\end{align}
where $s_{(\cdot)}=\text{sin}(\cdot),\, c_{(\cdot)}=\text{cos}(\cdot),\, t_{(\cdot)}=\text{tan}(\cdot)$, and $\lambda_{\Delta}=\lambda_{\Omega}+\lambda_{l}-\lambda_{t}$. To apply the Newton-Raphson method, we obtain the second derivative of $\lambda_{m}$ by differentiating Eq. \eqref{eq:dlam} for $t_b$ as follows:

\begin{align}
&\dfrac{\text{d}^{2}\lambda_{m}(t_b)}{\text{d}t_b^{2}}=-\frac{{c_{\beta}}\left(1+s^2_\phi\right)-2s_{\phi_t}s_\phi}{c_{\phi_t}c^3_{\phi}s_{\lambda_\Delta}}\dot\phi^2-\frac{c_{\beta}s_{\phi}-s_{\phi_t}}{{c_{\phi_t}c^3_{\phi}s_{\lambda_\Delta}}}\ddot\phi^2-\frac{(\dot{\lambda}_{\Omega}+\dot{\lambda}_{m}+\dot{\lambda}_{c})^2}{t_{\lambda_{\Delta}}}-\ddot{\lambda}_l-\ddot{\lambda}_c \label{eq:d2lam}
\end{align}

In Eqs. \eqref{eq:dlam}-\eqref{eq:d2lam}, the first and second derivatives of $\phi$, $\lambda_{l}$, and $\lambda_c$ are defined from Eqs. \eqref{eq:lat}, \eqref{eq:lon}, and \eqref{eq:lambda_c} as
\begin{equation}   
\begin{aligned}  
\dot{\phi} &= \dot{u} s_i \frac{c_u}{c_\phi},\quad \dot{\lambda}_l = \dot\Omega-\omega_{\oplus} +  \frac{\dot{u}c_i}{c^2_u+c^2_i s^2_u},\quad \dot{\lambda}_c=k_1 \\
\ddot{\phi} &= t_{\phi}(\dot{\phi}^{2}-\dot{u}^{2}), \quad \ddot{\lambda}_l = \frac{2 \dot{u}^2 t_u c_i s^{2}_i (1 + t^2_u)}{(1 + c^2_i t^2_u)^{2}},\quad \ddot{\lambda}_c=0
\end{aligned}
\end{equation}
Iteration proceeds by updating $t_b$ as $t_b=t_b-(\dot{\lambda}_m/\ddot{\lambda}_m)$ until $\dot{\lambda}_m$ reaches a predefined tolerance value. Since $\lambda_m$ has a shape similar to an ellipse, its extremes tend to occur near the center.  Therefore, selecting the midpoint, $(t_{L}+t_{U})/2$, as the initial estimate allows the Newton-Raphson method to converge to a highly accurate solution in just a small number of iterations (in most cases, three or four iterations are observed). Note that we should solve the root-finding problem separately for $\lambda_{m_U}$ and $\lambda_{m_L}$. The optimized solution provides the bounds of $\lambda_m$ ($\lambda_{m_L}$ and $\lambda_{m_U}$) and the corresponding time ($t_{m_L}$ and $t_{m_U}$). 

To enhance the accuracy of the model, we subsequently specify the four supplementary data points as
\begin{align}
\mathbf{t_s}=[t_{s_1}, t_{s_2}, t_{s_3}, t_{s_4}]
\end{align}
These points are distributed over the intervals formed by the six pre-determined data points, with priority given to relatively wider gaps so that sparsely sampled regions can be represented more accurately. Depending on the relative ordering of $t_{\Omega_U}$ and $t_{m_U}$, and on the relative distance from $t_{\Omega_U}$ to its neighboring boundary points, the placement strategy is classified into the following three cases:

\begin{itemize}
    \item \textbf{Case 1:} $t_{\Omega_U}<t_{m_U}$ and $|t_{\Omega_U}-t_{m_U}| > |t_{\Omega_U}-t_L|$
    \item \textbf{Case 2:} $t_{\Omega_U}<t_{m_U}$ and $|t_{\Omega_U}-t_{m_U}| \le |t_{\Omega_U}-t_L|$
    \item \textbf{Case 3:} $t_{\Omega_U}>t_{m_U}$
\end{itemize}
Accordingly, the supplementary points are determined as
\begin{align}
\mathbf{t}_{s} =
\begin{cases}
\left[
\dfrac{t_{\Omega_{U}} + t_{m_U}}{2},\;
\dfrac{t_{m_{U}} + t_{U}}{2},\;
\dfrac{t_{\Omega_{L}} + t_{m_L}}{2},\;
\dfrac{t_{m_{L}} + t_{L}}{2}
\right], 
& \text{for Case 1}, \\[1.2em]
\left[
\dfrac{t_{\Omega_{U}} + t_{L}}{2},\;
\dfrac{t_{m_U}+t_{U}}{2},\;
\dfrac{t_{\Omega_{L}} + t_{U}}{2},\;
\dfrac{t_{m_{L}} + t_{L}}{2}
\right], 
& \text{for Case 2}, \\[1.2em]
\left[
\dfrac{t_{L} + t_{m_U}}{2},\;
\dfrac{t_{\Omega_{U}} + t_{U}}{2},\;
\dfrac{t_U + t_{m_L}}{2},\;
\dfrac{t_{\Omega_{L}} + t_{L}}{2}
\right], 
& \text{for Case 3}.
\end{cases}
\label{eq:ts}
\end{align}

Consequently, we obtain ten time stamps composed of six representative points from the feasible range and four supplementary points. We collect them as
\begin{align}
\mathbf{t}_{\text{Data}}=[t_L, t_U, t_{\Omega_L},t_{\Omega_U},t_{m_L},t_{m_U},\mathbf{t}_s].
\end{align}
together with their corresponding boundary values $\lambda_m(\mathbf{t}_{\text{Data}})$.

\section{Explicit Expressions for Access Time}
This section derives explicit expressions for the visibility-boundary-crossing time, and consequently the access time, as functions of $\lambda_\Omega$. For a given set of orbital elements, the proposed procedure first identifies the feasible ranges of ELAN and MELAN and extracts a small number of representative boundary points, including selected intermediate points. These data points are then used to construct a piecewise ellipse representation of the MELAN boundary, $\lambda_m=f_i(t_b)$, where each segment is defined from adjacent pairs of data points and admits an analytic inverse. Using the defining relationship between ELAN and MELAN in Eq. \eqref{eq:MELAN}, 
the piecewise expression for $\lambda_m(t_b)$ is transformed into an explicit relation between $\lambda_\Omega$ and $t_b$. Because the transformed relation retains an analytically invertible quadratic form for a prescribed $\lambda_\Omega$, the inverse mapping $t_b(\lambda_\Omega)$ can be obtained explicitly for each selected segment and branch. The entry and exit times are then evaluated from this inverse relation, and the access time is subsequently computed.

\subsection{Alternative Formulation of $\lambda_m(t_b)$}
In this subsection, we derive an explicit expression of $\lambda_m(t_b)$ based on the proposed piecewise ellipse models. Using given data points $\mathbf{t}_{\text{Data}}$ and $\lambda_m(\mathbf{t}_{\text{Data}})$, the visible boundary is approximated by segmented elliptical curves, enabling an explicit formulation of $\lambda_m$ as a function of $t_b$. This formulation serves as the foundation for the subsequent transformation into $\lambda_\Omega(t_b)$ and its inverse mapping. 

Since the piecewise ellipse models are constructed in the two-dimensional planes (the $x-y$ plane for $t-\lambda_{m}$ and the $x-z$ plane for $t-\lambda_{\Omega}$), we substitute the variables with simplified notations as
\begin{align}
x := t_b, \qquad y := \lambda_m, \qquad z:=\lambda_{\Omega}.
\end{align}
Accordingly, all subsequent derivations are carried out in the $(x,y)$-plane. Then, the obtained data points also become
\begin{align}
\mathbf{x}_{\text{Data}} := \mathbf{t}_{\text{Data}}, 
\qquad
\mathbf{y}_{\text{Data}} := y(\mathbf{x}_{\text{Data}}).
\end{align}
To construct the piecewise ellipse model, $\mathbf{x}_{\text{Data}}$ is restructured so that the points are matched to the upper and lower visible boundaries, respectively. The data points for these two groups are defined as
\begin{equation}
\begin{aligned}
\mathbf{x}_{\text{Data}}^{+} &= [x_1, x_2, x_3, x_4, x_5, x_6], 
& \mathbf{x}_{\text{Data}}^{-} &= [x_7, x_8, x_9, x_{10}, x_{11}, x_{12}], \\
\mathbf{y}_{\text{Data}}^{+} &= y(\mathbf{x}_{\text{Data}}^{+}), 
& \mathbf{y}_{\text{Data}}^{-} &= y(\mathbf{x}_{\text{Data}}^{-}).
\end{aligned}
\end{equation}
where $x_1 = x_7 = t_L$ and $x_6 = x_{12} = t_U$, representing the visible range limits. The intermediate elements $x_2$--$x_5$ correspond to $\{t_{\Omega_U},\, t_{m_U},\, t_{s_1},\, t_{s_2}\}$ and are ordered in ascending $t$-value such that $x_2 < x_3 < x_4 < x_5$. Similarly, $x_8$--$x_{11}$ correspond to 
$\{t_{\Omega_L},\, t_{m_L},\, t_{s_3},\, t_{s_4}\}$ and are arranged in ascending order such that $x_8 < x_9 < x_{10} < x_{11}$. This construction ensures a consistent boundary representation in the $(x,y)$-plane for the subsequent piecewise ellipse formulation. The ordering of the data points in each group is determined by the location of $t_s$ in Eq. \eqref{eq:ts}.

Using this data set, we approximate $\lambda_m(t_b)$ with a ten-segment piecewise model whose each ellipse segment is defined by a pair of adjacent data points, and each group is composed of five such segments corresponding to consecutive portions of the curve as
\begin{align}
f_i(t_b,\lambda_m):=f_i(x,y)=
\displaystyle \frac{(x-x_{c})^2}{A_i^2}+\frac{y^2}{B_i^2}-1=0, \quad (i=1,2,...,10) \label{eq:ellipse}
\end{align}
where $i$ is the index number representing the $i$-th ellipse model shown in Fig. \ref{Fig:MELAN}, $x_{c}$ is the center of the $i$-th ellipse given as
\begin{align}
x_{c}=
\begin{cases}
t_{m_U}^+, \quad \text{for } i=1,2,3,4,5\\
t_{m_L}^-, \quad \text{for } i=6,7,8,9,10\\
\end{cases}
\end{align} Starting from the first data point $[x_1,0]$ (or $[x_7,0]$ for the lower boundary), the ellipse models are sequentially constructed in pairs of two points. Models 1 to 5 are defined by $\mathbf{x}_{\text{Data}}^+$, while Models 6 to 10 are defined by $\mathbf{x}_{\text{Data}}^-$. Moreover, $A_i$ and $B_i$ are the semi-major and semi-minor axes, respectively, calculated from the adjacent two points as
\begin{equation}
\begin{aligned}
A^2 &= 
\frac{(x_i-x_c)^2\,y_{i+1}^2-(x_{i+1}-x_c)^2\,y_i^2}{\,y_{i+1}^2-y_i^2\,}, \\[6pt]
B^2 &= 
\frac{(x_i-x_c)^2\,y_{i+1}^2-(x_{i+1}-x_c)^2\,y_i^2}{\,(x_i-x_c)^2-(x_{i+1}-x_c)^2\,}.
\end{aligned}
\end{equation}

\subsection{Analytic Expression for Visibility-Boundary-Crossing Time, $t_b(\lambda_{\Omega})$}
Eq. \eqref{eq:ellipse} provides the relationship between $t_b$ and $\lambda_m$. However, the goal of this subsection is to derive the explicit expression for $t_b$ as a function of $\lambda_{\Omega}$. Accordingly, we put Eqs. \eqref{eq:lambda_c} and \eqref{eq:MELAN} into Eq. \eqref{eq:ellipse} to explicitly determine the relationship between $t_b$ and $\lambda_{\Omega}$ as follows
\begin{align}
f_i(t_b,\lambda_{\Omega}):=&f_i(x,z)=\displaystyle \frac{(x-x_c)^2}{A_i^2}+\frac{(z-k_1x-k_2)^2}{B_i^2}-1 \nonumber\\
=&(B_i^2+A_i^2k_1^2)x^2 +2(-B_i^2x_c - A_i^2k_1z+A_i^2k_1k_2)x+B_i^2x_c^2+A_i^2z^2-2A_i^2k_2z+A_i^2k_2^2-A_i^2B_i^2 \label{eq:ellipse2}\\
=& \ 0 \nonumber
\end{align}
\begin{figure}[hbt!]
\centering
\includegraphics[width=0.7\columnwidth]{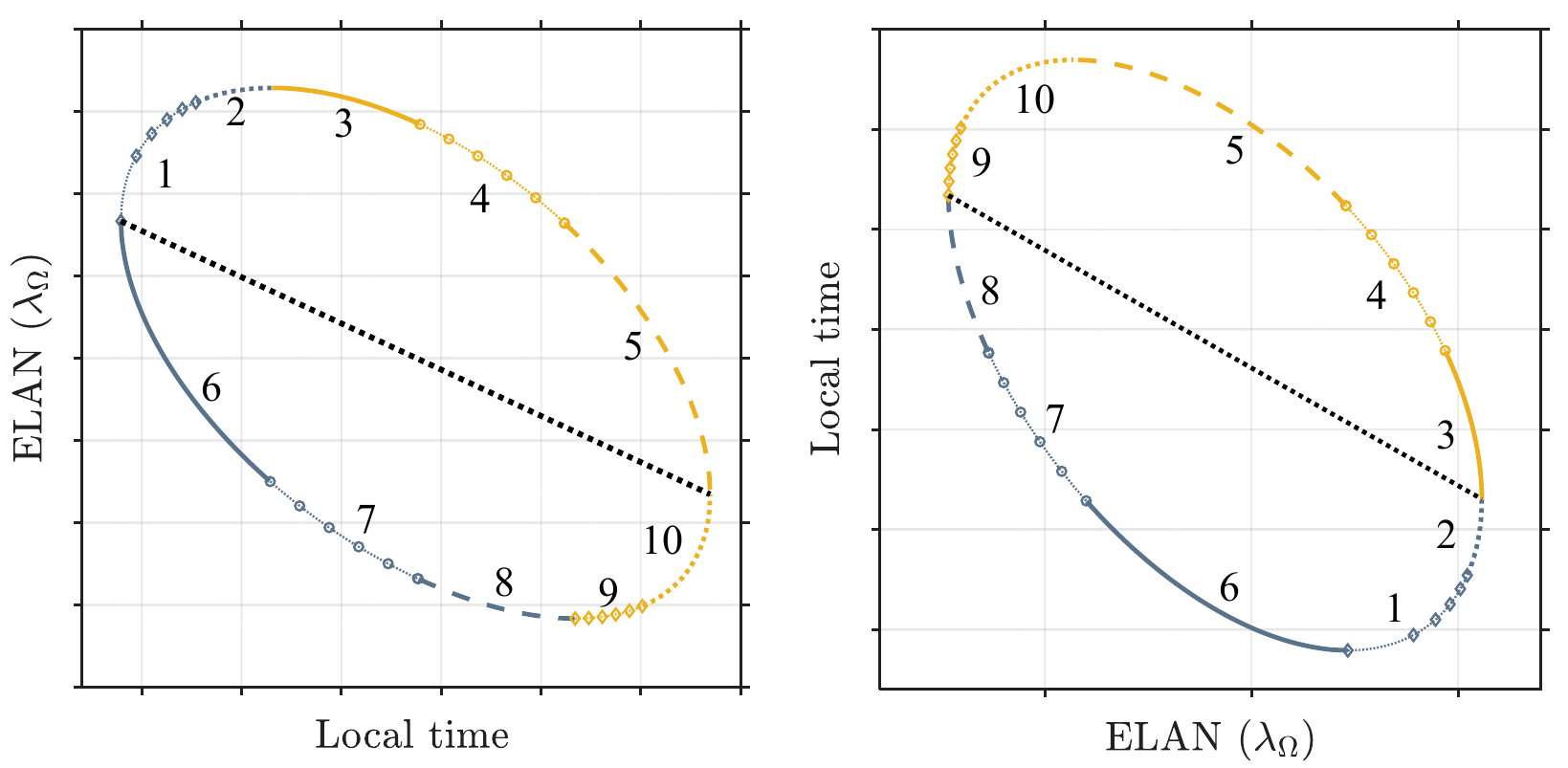}
\caption{Illustrative examples of $\lambda_{\Omega}(t_b)$ and $t_b(\lambda_{\Omega})$ obtained from the piecewise ellipse model (Case 2).}
\label{Fig_Visible_boundary_model}
\end{figure}

Eq. \eqref{eq:ellipse2} is a quadratic equation in $x(=t_b)$. Rearranging Eq. \eqref{eq:ellipse2} with respect to $x$ and applying the quadratic formula yields the explicit expression for $t_b$ as a function of $z(=\lambda_{\Omega})$:
\begin{align}
t_b =
\begin{cases}
t_{\mathrm{out}} = f_i^{+}(z) = \dfrac{-Q_i+\sqrt{Q_i^2-P_iR_i}}{P_i}, \\[0.8em]
t_{\mathrm{in}}  = f_i^{-}(z) = \dfrac{-Q_i-\sqrt{Q_i^2-P_iR_i}}{P_i},
\end{cases}
\label{eq:ta}
\end{align}
where
\begin{align}
P_i &= B_i^2 + A_i^2k_1^2, \nonumber\\
Q_i &= -B_i^2x_c - A_i^2k_1z + A_i^2k_1k_2, \\
R_i &= B_i^2x_c^2 + A_i^2z^2 - 2A_i^2k_2z + A_i^2k_2^2 - A_i^2B_i^2. \nonumber
\end{align}

The branch in Eq. \eqref{eq:ta} is selected according to the range of $z(=\lambda_{\Omega})$ associated with the $i$-th model. Figure \ref{Fig_Visible_boundary_model} presents illustrative examples of both $\lambda_{\Omega}(t_b)$ and its inverse function, $t_b(\lambda_{\Omega})$, generated using the piecewise ellipse model proposed in this paper for the graph shown in Fig. \ref{Fig_visible_boundary}. Specifically, the figure is obtained by plotting Eq. \eqref{eq:ellipse2} in the $t_b$--$\lambda_{\Omega}$ and $\lambda_{\Omega}$--$t_b$ coordinate systems, respectively. The model indices 1 through 10 in Fig. \ref{Fig_Visible_boundary_model} correspond to those shown in Fig. \ref{Fig:MELAN}. Depending on Cases 1--3, the ordering of the model indices (1--10) in the inverse representation, $t_b(\lambda_{\Omega})$, differs from that in the original representation, $\lambda_{\Omega}(t_b)$. Accordingly, the explicit expression for $t_b$ in Eq. \eqref{eq:ta} should be interpreted using both the selected branch and the corresponding model-index sequence for each case, given by
\begin{align}
\text{Case 1:}\qquad \label{eq:tindex}
&t_{\mathrm{out}}:\quad i \in \{10,\,5,\,4,\,3,\,2\},\quad t_{\mathrm{in}}:\quad i \in \{9,\,8,\,7,\,6,\,1\}, \nonumber \\[0.4em] 
\text{Case 2:}\qquad
&t_{\mathrm{out}}:\quad i \in \{9,\,10,\,5,\,4,\,3\}, \quad t_{\mathrm{in}}:\quad i \in \{8,\,7,\,6,\,1,\,2\}, \\[0.4em] 
\text{Case 3:}\qquad
&t_{\mathrm{out}}:\quad i \in \{8,\,9,\,10,\,5,\,4\}, \quad t_{\mathrm{in}}:\quad i \in \{7,\,6,\,1,\,2,\,3\}. \nonumber 
\end{align}

\subsection{Coverage Analysis}
With given initial orbital elements and $\lambda_{\Omega}(t_0)(=\Omega_0-\theta_{G_0})$, we can generate the ELAN profile over the entire analysis horizon using Eq.~\eqref{eq:ELAN} as
\begin{align}
\boldsymbol{\lambda}_{\Omega} = \big[\lambda_{\Omega}^{(1)}, \lambda_{\Omega}^{(2)}, \ldots, \lambda_{\Omega}^{(n_f)}\big],
\end{align}
For each ELAN $\lambda_{\Omega}^{(n)}$, the proposed procedure first performs a feasibility check to determine whether the corresponding geometry lies within the visible range. If $\lambda_{\Omega}^{(n)}$ is infeasible (i.e., outside the visible range), the access time for that nodal period is set to zero. Otherwise, $\lambda_{\Omega}^{(n)}$, determined to be feasible, is substituted into Eqs.~\eqref{eq:ta}--\eqref{eq:tindex} to compute the entry and exit times, $t_{\mathrm{in}}^{(n)}$ and $t_{\mathrm{out}}^{(n)}$, respectively. The access time within the $n$-th nodal period is then evaluated as
\begin{align}
t_a^{(n)} = t_{\mathrm{out}}^{(n)} - t_{\mathrm{in}}^{(n)}, \qquad n=1,2,\ldots,n_f.
\label{Eq: ta_n}
\end{align}
For the last nodal period, if the local final time, \(t_f^{(n_f)}\), is smaller than \(t_{\mathrm{out}}^{(n_f)}\), the access time is determined by the following two cases:
\begin{align}
t_a^{(n_f)} =
\begin{cases}
t_f^{(n_f)} - t_{\mathrm{in}}^{(n_f)}, & t_{\mathrm{in}}^{(n_f)} < t_f^{(n_f)} < t_{\mathrm{out}}^{(n_f)}, \\[4pt]
0, & t_f^{(n_f)} < t_{\mathrm{in}}^{(n_f)}.
\end{cases}
\end{align}
Repeating this process for the ELAN profile yields the access-time sequence
$\big[t_a^{(1)}, t_a^{(2)}, \ldots, t_a^{(n_f)}\big]$.
Finally, the total access time over the entire analysis horizon is obtained by summing the access times across all nodal periods:
\begin{align}
t_{a,\text{tot}} = \sum_{n=1}^{n_f} t_a^{(n)}.
\end{align}

To evaluate the revisit performance, the local entry and exit times are converted into absolute times over the entire time horizon as
\begin{align}
\tilde{t}_{\mathrm{in}}^{(n)} &= (n-1)T_n + t_{\mathrm{in}}^{(n)}, \\
\tilde{t}_{\mathrm{out}}^{(n)} &= (n-1)T_n + t_{\mathrm{out}}^{(n)},
\qquad n=1,2,\ldots,n_f.
\end{align}
Only the feasible access windows satisfying $t_a^{(n)} > 0$ are retained for revisit-time evaluation. Let the ordered sequence of feasible access windows be
\begin{align}
\left\{
\left(\tilde{t}_{\mathrm{in}}^{(1)}, \tilde{t}_{\mathrm{out}}^{(1)}\right),
\left(\tilde{t}_{\mathrm{in}}^{(2)}, \tilde{t}_{\mathrm{out}}^{(2)}\right),
\ldots,
\left(\tilde{t}_{\mathrm{in}}^{(N_w)}, \tilde{t}_{\mathrm{out}}^{(N_w)}\right)
\right\},
\end{align}
where $N_w$ denotes the number of feasible access windows, equivalent to the number of visible ELANs. The revisit time between two consecutive feasible access intervals is then defined as
\begin{align}
t_r^{(k)} = \tilde{t}_{\mathrm{in}}^{(k+1)} - \tilde{t}_{\mathrm{out}}^{(k)},
\qquad k = 1,2,\ldots,N_w-1.
\end{align}
Accordingly, the maximum revisit time over the entire analysis period is obtained as
\begin{align}
t_{r,\max} = \max_{k=1,\ldots,N_w-1}
\left(
\tilde{t}_{\mathrm{in}}^{(k+1)} - \tilde{t}_{\mathrm{out}}^{(k)}
\right).
\label{Eq: maxrt}
\end{align}

\section{Case Studies: Validation and Application to Coverage Analysis}
This section presents case studies to demonstrate both the accuracy and practical usefulness of the proposed approach. First, the proposed explicit visibility-boundary-crossing time expression with respect to ELAN is validated against the simulation results under various orbital and observation conditions. This validation quantitatively examines how accurately the proposed model approximates the relationship between ELAN and access time over the visible range of ELAN. Second, the proposed approach is applied to coverage analysis to assess its practical performance in terms of both prediction accuracy and computational efficiency. 

\subsection{Validation of the proposed method}
To validate the proposed approach, we compare the access time obtained from the proposed approach with that obtained from a simplified vectorization-based method used in \cite{lee2020satellite} over the entire visible ELAN interval ($[\lambda_{\Omega_L}, \lambda_{\Omega_U}]$). In the reference method, the orbital dynamics are modeled using mean orbital elements that incorporate the J2 perturbation, and the satellite longitude and latitude are obtained analytically as functions of time. The access time is then determined by discretizing time and evaluating the visibility condition at each time step.  Specifically, for each ELAN sample within the visible range, the ground-track profile in Eqs.~\eqref{eq:u}--\eqref{eq:lon} is generated in a fully vectorized manner over time, and the visibility function in Eq.~\eqref{eq:NewV} is evaluated at each time step. The access time is computed by identifying the time intervals during which the visibility condition is satisfied. 

The visibility geometry is determined by the satellite altitude, orbital inclination, target latitude, and sensor half-angle. Since these parameters jointly affect the size and shape of the visible region, the accuracy of the proposed method is examined over a wide range of conditions. For each orbital configuration, the visible ELAN bounds are first determined, and the interval is discretized into sufficiently dense ELAN samples. The access time at each ELAN sample is then computed using both the vectorized method and the proposed approach, and the errors are evaluated over all ELAN samples within the visible range.

We consider two types of error metrics: 1) the mean pointwise relative error, obtained by computing the relative error at each ELAN sample and then averaging it over the entire visible ELAN interval and 2) the integrated relative error, obtained by first summing the access times predicted by the proposed model over all visible ELAN samples, then summing the corresponding exact access times, and finally computing the relative error between these two aggregated values. 
\begin{align}
\varepsilon_{\mathrm{1}}=\frac{1}{N_{\mathrm{sp}}}\sum_{s=1}^{N_{\mathrm{sp}}}\left|\frac{\hat{t}_{a,s}-t_{a,s}}{t_{a,s}}\right|,
\qquad
\varepsilon_{\mathrm{2}}=\left|\frac{\sum_{s=1}^{N_{\mathrm{sp}}}\hat{t}_{a,s}-\sum_{s=1}^{N_{\mathrm{sp}}}t_{a,s}}{\sum_{s=1}^{N_{\mathrm{sp}}}t_{a,s}}\right|
\
\end{align}
where \(\hat{t}_{a,s}\) and \(t_{a,s}\) denote the true and predicted access times at the \(s\)-th visible ELAN sample, respectively, and \(N_{\mathrm{sp}}\) is the total number of visible ELAN samples. In this study, \(N_{\mathrm{sp}}\) is set to 50 and the time step is set to 0.01 seconds. The second metric is introduced since the relative error of the access time varies significantly depending on the ELAN value. Near the lower or upper visible ELAN bounds, the true access time may be only a few seconds, and thus even a small absolute discrepancy can produce an excessively large pointwise relative error. In contrast, the integrated relative error better reflects the overall agreement in the access-time distribution across the full visible ELAN range.
\begin{table*}[t]
\centering
\caption{Validation results of the proposed approach ($\varepsilon_1/\varepsilon_2$ in \%).}
\label{tab:validation}

\begingroup
\fontsize{10}{12}\selectfont
\renewcommand{\arraystretch}{1.15}
\setlength{\tabcolsep}{0pt}

\newcolumntype{C}[1]{>{\centering\arraybackslash}m{#1}}

\begin{tabular}{
C{1.55cm}
C{1.52cm} C{1.52cm} C{1.52cm} C{1.52cm} C{1.52cm}
C{1.52cm} C{1.52cm} C{1.52cm} C{1.52cm} C{1.52cm}
}
\toprule

Altitude & \multicolumn{5}{c}{Inclination [deg]} & \multicolumn{5}{c}{Inclination [deg]} \\
{[km]} & 20 & 40 & 60 & 80 & 100 & 20 & 40 & 60 & 80 & 100 \\
\midrule
& \multicolumn{5}{c}{$\phi_t=0^\circ$}
& \multicolumn{5}{c}{$\phi_t=20^\circ$} \\
\cmidrule(lr){2-6}\cmidrule(lr){7-11}
300
& 0.02/0.01 & 0.03/0.00 & 0.03/0.00 & 0.03/0.00 & 0.03/0.00
& -- & 0.04/0.01 & 0.03/0.01 & 0.03/0.00 & 0.04/0.00 \\

700
& 0.01/0.00 & 0.01/0.00 & 0.01/0.00 & 0.01/0.00 & 0.01/0.00
& -- & 0.04/0.00 & 0.02/0.00 & 0.02/0.00 & 0.02/0.00 \\

1,000
& 0.01/0.01 & 0.01/0.00 & 0.01/0.00 & 0.01/0.00 & 0.01/0.00
& -- & 0.06/0.00 & 0.02/0.00 & 0.02/0.00 & 0.01/0.00 \\

2,000
& 0.04/0.03 & 0.00/0.00 & 0.00/0.00 & 0.01/0.00 & 0.00/0.00
& -- & 0.11/0.01 & 0.04/0.00 & 0.02/0.00 & 0.02/0.00 \\

5,000
& 0.24/0.22 & 0.01/0.01 & 0.02/0.02 & 0.04/0.03 & 0.03/0.03
& -- & 0.30/0.10 & 0.11/0.01 & 0.07/0.03 & 0.06/0.02 \\

10,000
& 1.90/1.61 & 0.08/0.06 & 0.09/0.08 & 0.14/0.14 & 0.11/0.10
& -- & 1.06/0.81 & 0.30/0.03 & 0.19/0.13 & 0.14/0.10 \\

\midrule

& \multicolumn{5}{c}{$\phi_t=40^\circ$}
& \multicolumn{5}{c}{$\phi_t=60^\circ$} \\
\cmidrule(lr){2-6}\cmidrule(lr){7-11}
300
& -- & -- & 0.03/0.01 & 0.03/0.01 & 0.03/0.00
& -- & -- & -- & 0.04/0.00 & 0.04/0.00 \\

700
& -- & -- & 0.05/0.00 & 0.02/0.00 & 0.02/0.00
& -- & -- & -- & 0.04/0.00 & 0.03/0.01 \\

1000
& -- & -- & 0.07/0.00 & 0.03/0.00 & 0.03/0.00
& -- & -- & -- & 0.07/0.00 & 0.06/0.00 \\

2,000
& -- & -- & 0.13/0.01 & 0.06/0.00 & 0.05/0.00
& -- & -- & -- & 0.15/0.01 & 0.12/0.00 \\

5,000
& -- & -- & 0.35/0.11 & 0.18/0.02 & 0.13/0.02
& -- & -- & -- & 0.44/0.05 & 0.31/0.03 \\

10,000
& -- & -- & 0.73/0.68 & 0.47/0.06 & 0.28/0.07
& -- & -- & -- & 1.25/0.46 & 0.79/0.18 \\
\cmidrule(lr){1-11}
\label{Table: Validation}
\end{tabular}
\endgroup
\end{table*}

Table \ref{Table: Validation} summarizes the validation results of the proposed approach. In particular, the results are organized primarily with respect to the target latitude. For each target latitude, the errors are represented for multiple combinations of altitude and inclination. The rows correspond to altitudes of 300, 700, 1,000, 2,000, 5,000, and 10,000 km, and the columns correspond to inclinations of 20, 40, 60, 80, and 100 degrees. Four panels are provided for target latitudes of $0^\circ$, $20^\circ$, $40^\circ$, and $60^\circ$, respectively. We fix the sensor half-angle at $10^\circ$ as the variation of the visible region induced by the sensor half-angle is qualitatively similar to that induced by the orbital altitude; therefore, the influence of visibility-region size can be sufficiently represented by varying the altitude while keeping the sensor half-angle fixed.

Dashes in the tables correspond to physically infeasible cases in which the target is invisible under the given orbital configuration. The model accuracy is influenced by the orbital altitude, inclination, and target latitude, which jointly determine the structure of the visibility geometry. More specifically, the error exhibits a clear dependence on the duration for which the satellite remains within the visible region. As the orbital altitude increases, the satellite passes through the visible region for a longer duration, leading to a gradual increase in the error; accordingly, the largest errors are observed at the highest altitude. In contrast, as the inclination increases, the duration of passage through the visible region becomes shorter, resulting in a reduction in the error. On the other hand, the target latitude is found to have a relatively minor influence on the overall error trend.

\begin{figure}[hbt!]
\centering
\includegraphics[width=0.48\columnwidth]{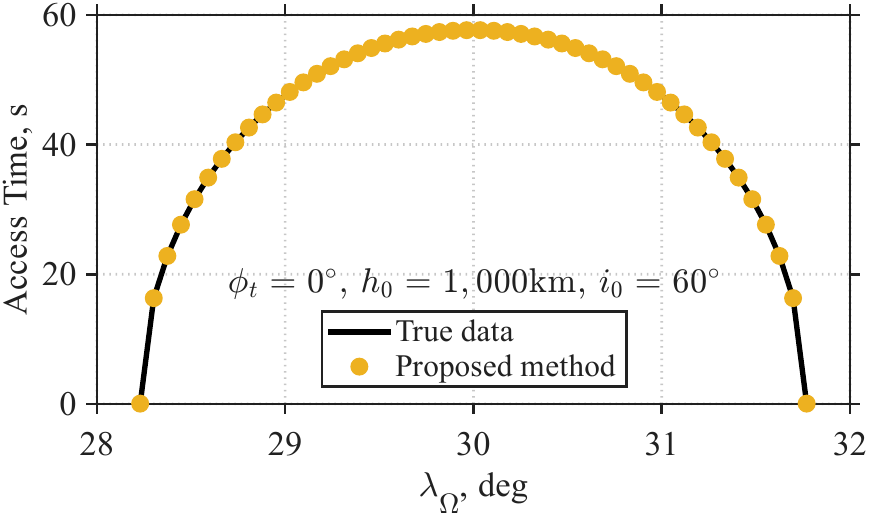}
\includegraphics[width=0.48\columnwidth]{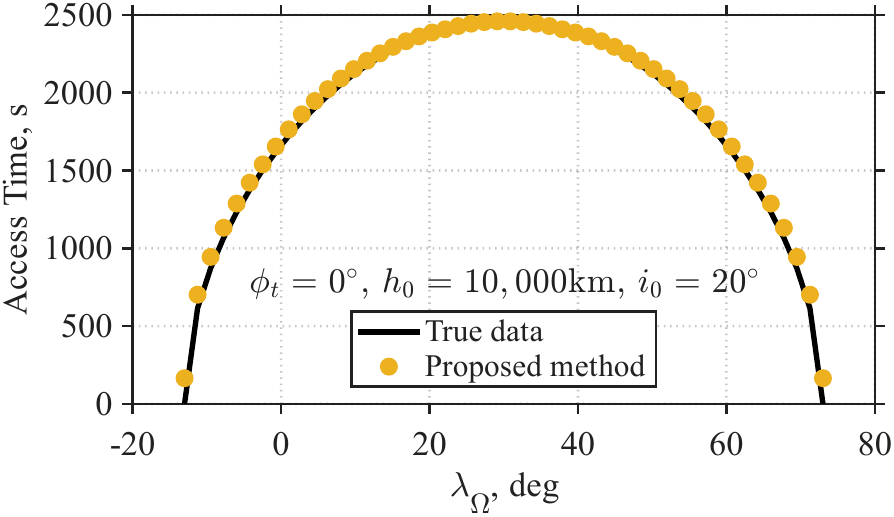}
\caption{Comparison results of the proposed approach}
\label{Fig_results1}
\end{figure}
To provide a more intuitive validation of the proposed approach, we present access-time profile versus ELAN for several selected combinations of target latitude, orbital altitude, and inclination within the visible ELAN bounds in Fig. \ref{Fig_results1}. For each case, the predicted access time from the proposed approach is compared directly with the true values. It is observed that the predicted points closely match the true values. 

The results show that the proposed model is highly accurate over a broad range of circular-orbit configurations. In most cases, particularly in low Earth orbit conditions, both the mean pointwise and integrated relative errors remain below $0.1\%$, indicating that the proposed method provides an almost exact approximation of the access time. The largest errors occur in high-altitude, low-inclination cases, where the target remains inside the visible region for a comparatively long duration and the fitted curve becomes flatter near its extrema. Even in those cases, however, the error remains sufficiently small for preliminary design and repeated trade studies.

 \subsection{Application to Coverage Analysis}
To further demonstrate the applicability of the proposed method, it is extended from the access-time analysis within a single nodal period to the coverage analysis over the entire analysis time for a single satellite. Specifically, the total target access time, maximum revisit time, and mean revisit time are evaluated over the entire period of interest. In the proposed approach, the ELAN profile is first generated over the entire analysis interval. Among these ELAN values, only those lying within the visible ELAN bounds are regarded as feasible, and the explicit expressions are applied only to those feasible cases to compute the corresponding entry/exit times. The overall coverage metrics over the entire analysis time are then obtained using Eqs.~\eqref{Eq: ta_n}--\eqref{Eq: maxrt}, including the total access time and the maximum revisit time. These results are finally compared with those obtained from the vectorized simulation-based method.

For the vectorized simulation-based method, two time-step values are considered: 0.1 s and 1 s. Among these, the result obtained with the finer time step of 0.1 s is regarded as a reference solution.

We perform a coverage analysis for a single satellite under three representative orbital configurations, denoted as Case I (representative case), Case II (best-case scenario), and Case III (worst-case scenario). Case I corresponds to a medium-altitude, medium-inclination orbit (5,000 km, 60 degrees). Case II considers a low-altitude, high-inclination orbit (300 km, 80 degrees), where the access time within the visible region is relatively short. In contrast, Case III considers a high-altitude, low-inclination orbit (10,000 km, 20 degrees), where the satellite traverses the visible region for a longer duration.

For Cases I and II, the target is located at a latitude of 30 degrees and, for Case III, it is set to 0 degrees. For all cases, a longitude of 30 degrees and the sensor half-angle of 10 degrees are considered. To capture both short-term and long-term coverage characteristics, the analysis periods are set to two weeks and one year, respectively.

\begin{figure}[hbt!]
\centering
\includegraphics[width=0.48\columnwidth]{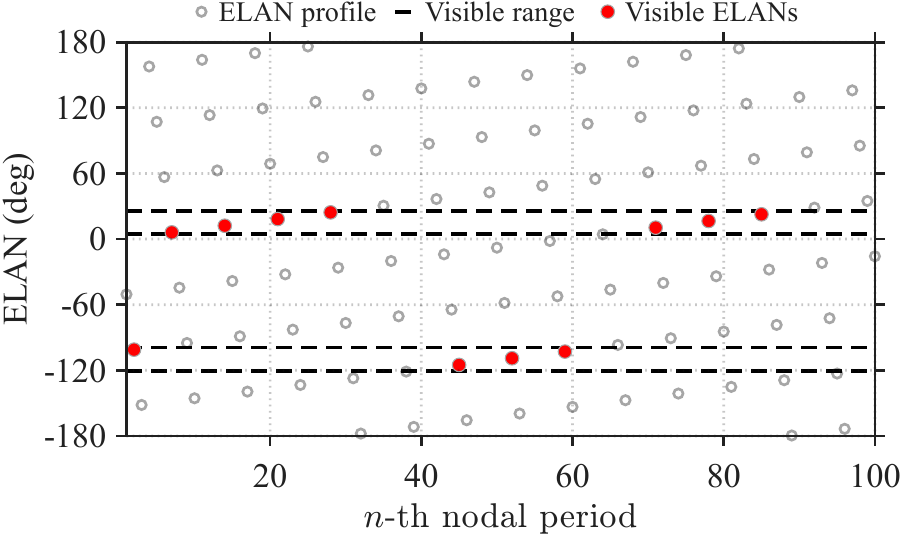}
\includegraphics[width=0.48\columnwidth]{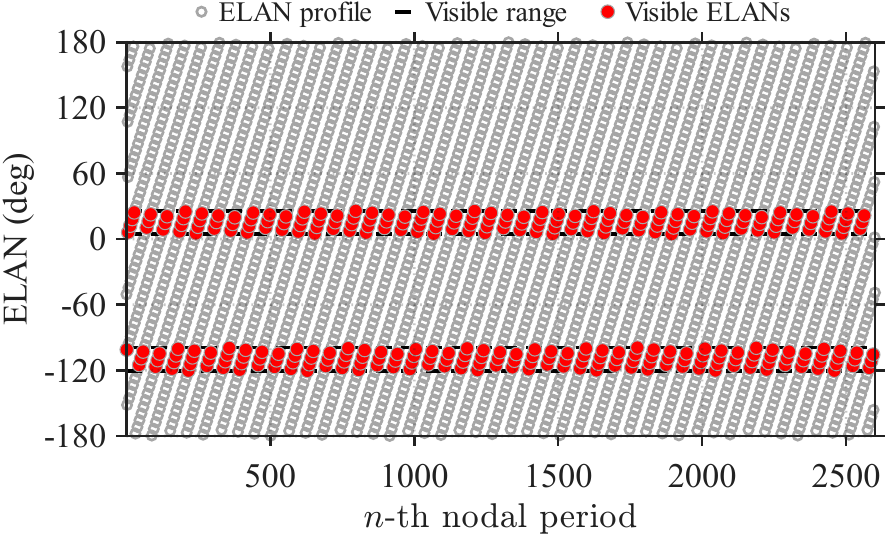}
\caption{ELAN profile versus $n$-th nodal period: 2 weeks (left) and 1 year (right) analysis (Case I)}
\label{Fig_result_ELAN_profile}
\end{figure}

\begin{figure}[hbt!]
\centering
\includegraphics[width=0.48\columnwidth]{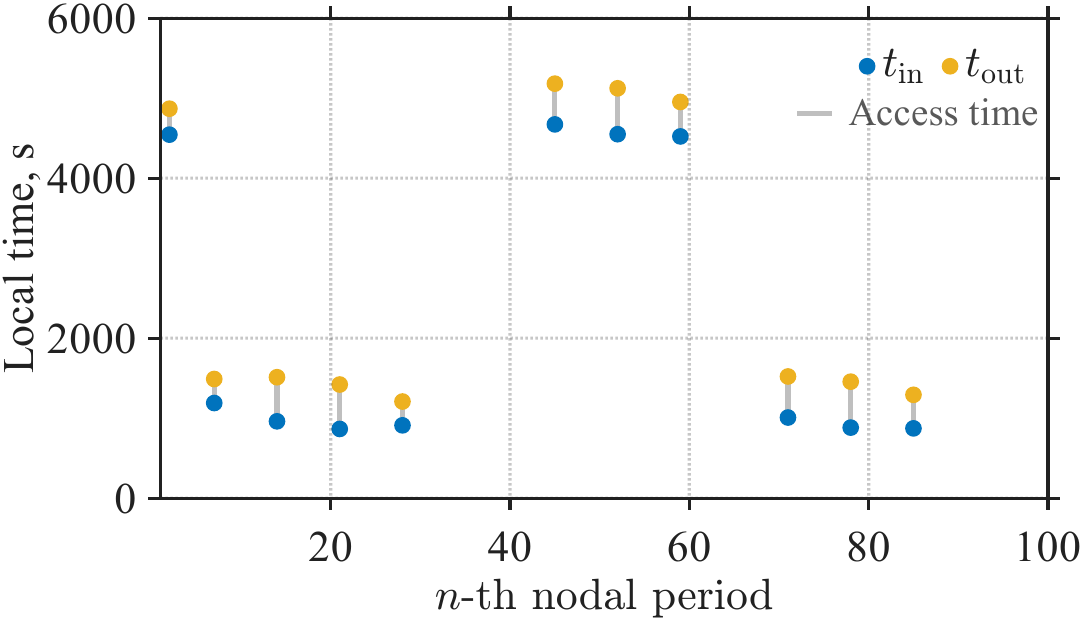}
\caption{Visibility-boundary-crossing time profile versus $n$-th nodal period: 2 weeks analysis (Case I)}
\label{Fig_result_access_epoch}
\end{figure}

For each case, the proposed approach first generates the ELAN profiles over the entire time span and identifies the subset of visible ELANs. Fig. \ref{Fig_result_ELAN_profile} shows the ELAN profile and visible ELANs for Case I. In Fig. \ref{Fig_result_ELAN_profile}, two sets of visible range exist as the orbit could enter the visible region along the ascending and descending tracks. After selecting the visible ELAN values, the visibility-boundary-crossing time corresponding to those values is calculated from the proposed method. Fig. \ref{Fig_result_access_epoch} presents the visibility-boundary-crossing time profile associated with the visible ELAN values for a two-week analysis for Case I, distinguishing between the entry and exit times, denoted by $t_{\text{in}}$ and $t_{\text{out}}$, respectively. From the profile, we can compute the total access time and maximum revisit time during the entire time span.

\begin{table}[!t]
\centering
\caption{Comparison of coverage analysis results (2-week analysis)}
\begin{tabular}{c c c c c c c}
\hline
\multirow{2}{*}{Case} & \multirow{2}{*}{Metric} 
& \multicolumn{2}{c}{Vectorized method} 
& \multirow{2}{*}{\makecell{Proposed\\Method}}
& \multirow{2}{*}{Abs. error} 
& \multirow{2}{*}{Rel. error (\%)} \\
\cline{3-4}
& & $\Delta t = 0.1$ s & $\Delta t = 1$ s & & & \\
\hline

\multirow{4}{*}{\begin{tabular}{c}Case I\\(5,000 km, 60 deg)\end{tabular}}
& Total access time [s]   & 5,051.7 & 5,050.0 & 5,051.1 & 0.6 & 0.0 \\
& Mean revisit time [hr]  & 27.6 & 27.6 & 27.6 & 0.0 & 0.0 \\
& Max revisit time [hr]   & 58.0 & 58.0 & 58.0 & 0.0 & 0.0 \\
& Computation time [ms]   & 330 & 36 & 0.3 & -- & -- \\
\hline

\multirow{4}{*}{\begin{tabular}{c}Case II\\(300 km, 80 deg)\end{tabular}}
& Total access time [s]   & 19.7 & 20.0 & 19.8 & 0.1 & 0.5 \\
& Mean revisit time [hr]  & 203.4 & 203.4 & 203.4 & 0.0 & 0.0 \\
& Max revisit time [hr]   & 203.4 & 203.4 & 203.4 & 0.0 & 0.0 \\
& Computation time [ms]   & 315 & 37 & 0.2 & -- & -- \\
\hline

\multirow{4}{*}{\begin{tabular}{c}Case III\\(10,000 km, 20 deg)\end{tabular}}
& Total access time [s]   & 53,273 & 53,271 & 52,556 & 717 & 1.3 \\
& Mean revisit time [hr]  & 11.2 & 11.2 & 11.4 & 0.2 & 1.8 \\
& Max revisit time [hr]   & 14.9 & 14.9 & 14.8 & 0.1 & 0.6 \\
& Computation time [ms]   & 332 & 34 & 0.3 & -- & -- \\
\hline
\end{tabular}
\label{Table:comparison1}

\vspace{20pt} 

\centering
\caption{Comparison of coverage analysis results (1-year analysis)}
\begin{tabular}{c c c c c c c}
\hline
\multirow{2}{*}{Case} & \multirow{2}{*}{Metric} 
& \multicolumn{2}{c}{Vectorized method} 
& \multirow{2}{*}{\makecell{Proposed\\Method}}
& \multirow{2}{*}{Abs. error} 
& \multirow{2}{*}{Rel. error (\%)} \\
\cline{3-4}
& & $\Delta t = 0.1$ s & $\Delta t = 1$ s & & \\
\hline

\multirow{4}{*}{\begin{tabular}{c}Case I\\(5,000 km, 60 deg)\end{tabular}}
& Total access time [s]   & 140,728 & 140,724 & 140,727 & 1 & 0.0 \\
& Mean revisit time [hr]  & 28.0 & 28.0 & 28.0 & 0.0 & 0.0 \\
& Max revisit time [hr]   & 58.0 & 58.0 & 58.0 & 0.0 & 0.0 \\
& Computation time [ms]   & 7,900 & 826 & 1.2 & -- & -- \\
\hline

\multirow{4}{*}{\begin{tabular}{c}Case II\\(300 km, 80 deg)\end{tabular}}
& Total access time [s]   & 411 & 414 & 412 & 1 & 0.2 \\
& Mean revisit time [hr]  & 246.9 & 246.9 & 246.9 & 0.0 & 0.0 \\
& Max revisit time [hr]   & 559.4 & 559.4 & 559.4 & 0.0 & 0.0 \\
& Computation time [ms]   & 8,930 & 912 & 0.5 & -- & -- \\
\hline

\multirow{4}{*}{\begin{tabular}{c}Case III\\(10,000 km, 20 deg)\end{tabular}}
& Total access time [s]   & 1,356,336 & 1,356,332 & 1,377,892 & 21,556 & 1.6 \\
& Mean revisit time [hr]  & 11.6 & 11.6 & 11.6 & 0.0 & 0.0 \\
& Max revisit time [hr]   & 22.3 & 22.3 & 22.3 & 0.0 & 0.0 \\
& Computation time [ms]   & 8,740 & 830 & 3 & -- & -- \\
\hline

\end{tabular}
\label{Table:comparison2}
\end{table}
\FloatBarrier

Tables~\ref{Table:comparison1} and \ref{Table:comparison2} summarize the results of the coverage analysis for the two-week and one-year periods, respectively, including comparisons with the simplified-vectorized method. The results obtained with the finer time step of 0.1 s are treated as the reference solution, and all errors are computed with respect to this reference. As expected, the vectorized simulation-based method becomes more accurate as the time-step size decreases; however, this improvement comes at the cost of increased computational time. The computation time is reported as the average value obtained from ten independent runs to account for variability in execution time in MATLAB R2025b on a desktop computer equipped with an AMD Ryzen 9 5950X 16-core processor and 64 GB RAM.

Overall, the proposed approach provides highly accurate results with substantially reduced computational cost compared to the vectorized simulation-based method, demonstrating its effectiveness for efficient coverage analysis. When applied to coverage analysis over extended time horizons, including long-term scenarios, the proposed method enables rapid computation of key performance metrics such as the total access time, mean revisit time, and maximum revisit time without compromising accuracy. Furthermore, the results demonstrate that the proposed approach remains applicable across a wide range of orbital element combinations, highlighting its robustness and flexibility. This capability makes it particularly well suited for orbit design applications, where it can significantly reduce the computational burden associated with repeated coverage evaluations.

\section{Conclusion}
This paper proposes a semi-analytic method for efficient coverage analysis of a satellite using piecewise ellipse models. The main idea is to exploit the invariance of the ground-track shape over one nodal period and to parameterize target visibility by the epoch longitude of the ascending node (ELAN). The feasible ELAN intervals were first identified, after which the relationships between ELAN and the entry/exit times were approximated with piecewise ellipse models and analytically inverted to obtain explicit visibility-boundary-crossing time expressions.
Once the visible ELAN intervals and the corresponding model coefficients are available, the proposed method can evaluate access time and revisit metrics over analysis horizons without repeated time-domain visibility searches. Numerical validation showed that the method achieves high accuracy across a broad range of circular-orbit configurations, with sub-0.1\% error in most LEO/MEO cases, while reducing computation time by several orders of magnitude relative to a vectorized time-stepped baseline. Extending the framework to elliptic orbits, higher-fidelity perturbation models, and constellations would be valuable directions for future research.

\section*{Acknowledgments}
This work was supported by the National Research Foundation of Korea(NRF) grant funded by the Korea government(MSIT) (RS-2025-24683783) and by the Korea Aerospace Administration (KASA) grant funded by the Korea government (RS-2026-25522528).

\bibliography{sample}

\end{document}